\documentclass[pdflatex,sn-mathphys-num]{sn-jnl}
\usepackage{amsmath,amssymb,amsthm,mathtools,booktabs,array}
\usepackage[T1]{fontenc}
\usepackage{lmodern}
\usepackage{microtype}
\usepackage{enumitem}
\hypersetup{hidelinks}
\usepackage{xcolor}

\theoremstyle{thmstyleone}
\newtheorem{theorem}{Theorem}[section]
\newtheorem{proposition}[theorem]{Proposition}
\newtheorem{lemma}[theorem]{Lemma}
\newtheorem{corollary}[theorem]{Corollary}
\newtheorem{conjecture}[theorem]{Conjecture}
\theoremstyle{thmstylethree}
\newtheorem{definition}[theorem]{Definition}
\theoremstyle{thmstyletwo}
\newtheorem{remark}[theorem]{Remark}

\newcommand{\F}{\mathbb F}
\newcommand{\C}{\mathbb C}
\newcommand{\A}{\mathbb A}

\newcommand{\Hom}{\operatorname{Hom}}
\newcommand{\End}{\operatorname{End}}
\newcommand{\Bil}{\operatorname{Bil}}
\newcommand{\Ind}{\operatorname{Ind}}
\newcommand{\Tr}{\operatorname{Tr}}

\newcommand{\im}{\operatorname{im}}
\newcommand{\Span}{\operatorname{span}}

\begin{document}

\title[Secant dimensions and bilinear widths]{Secant dimensions and bilinear decomposition widths for finite Heisenberg groups}
\author*[1]{\sur{Nguyen} \fnm{Bich Van}}\email{nbvan@vnu.edu.vn}
\author[1]{\sur{Nguyen} \fnm{Trung Kien}}\email{25022890@vnu.edu.vn}
\affil[1]{\orgdiv{Institute for Artificial Intelligence}, \orgname{VNU University of Engineering and Technology}, \orgaddress{\city{Hanoi}, \country{Vietnam}}}

\abstract{Let $G=UT_3(\F_q)$ be the finite Heisenberg group. We study the image closures and secant varieties of a family of equivariant bilinear maps between representations induced from central characters. Evaluation at a basepoint identifies the ambient space with $M_{q^2}(\C)$. A single channel is represented by $A^T\operatorname{diag}(d)B$, with $A$ and $B$ in prescribed twisted-convolution algebras. We prove that its affine image closure has dimension $3q^2-2$ for every central-character pair and that its affine second-secant closure has the expected dimension for every prime power $q$. For $q\ge3$, the proof constructs two differential images with zero intersection. We also prove a subgroup-lifting theorem for channel parameterizations on finite abelian groups with bicharacter coefficients: a parameter point with differential rank $r(3|H|-2)$ on a subgroup $H$ extends by zero to a parameter point with differential rank $r(3|K|-2)$ on the ambient group $K$, provided that every parameter vector is nonzero. Exact Jacobian certificates determine all secant dimensions for $q=2,3,4,5$, with generic widths $2,4,6,9$. The lifting theorem yields specified higher-secant dimensions on infinite families of extension fields. Attainment of the dimension lower bound for arbitrary $q$ remains a conjecture.}

\maketitle

\section{Introduction}

Let $q$ be a prime power and let $G=UT_3(\F_q)$ be the group of upper unitriangular $3\times3$ matrices over $\F_q$. Its center $Z$ is isomorphic to the additive group of $\F_q$. Fix a nontrivial additive character $\psi$ of $\F_q$, and put $\psi_t(z)=\psi(tz)$ for $t\in\F_q$. We consider the induced representations
\[
\mathcal H_t=\Ind_Z^G\psi_t.
\]
Each $\mathcal H_t$ has dimension $q^2$. Pointwise multiplication defines an equivariant bilinear map from $\mathcal H_s\times\mathcal H_t$ to $\mathcal H_{s+t}$. Thus the central characters determine the possible input and output types of a bilinear map.

Equivariant linear maps between these representations can be described by twisted convolutions. This is a finite version of the twisted-convolution calculus for Heisenberg representations; see Folland~\citep{Folland1989}. The relevant representation theory is part of finite Stone--von Neumann theory \citep{GurevichHadani2009,GurevichHadani2012}. In the preprint of Bich Van Nguyen and Nguyen Cao Manh Thang (2025; \href{https://arxiv.org/abs/2507.08397}{arXiv:2507.08397}), explicit decompositions of linear equivariant feature spaces and filters are studied for $UT_3(\F_3)$. Here we consider sums of bilinear maps of the form
\[
\Gamma_{A,B,C}(f,g)=C[(Af)(Bg)],
\]
where $A$, $B$, and $C$ are equivariant linear operators on $\mathcal H_s$, $\mathcal H_t$, and $\mathcal H_{s+t}$, respectively. We call such a map a fusion channel. The fusion width of an equivariant bilinear map is the least number of channels whose sum is that map.

Evaluation at a basepoint identifies the space of these equivariant bilinear maps with $M_{q^2}(\C)$. In these coordinates, a channel has the form $A^T\operatorname{diag}(d)B$, where $d$ is arbitrary and $A,B$ belong to the corresponding commutant algebras. Let $X_{s,t}$ be the affine closure of the one-channel image. The closure of the sum of $r$ channels is the affine cone over $\sigma_r(\mathbb P X_{s,t})$. We determine its dimension at order two for every prime power and every pair of central characters. We also prove a lifting theorem for parameterizations defined by bicharacters on finite abelian groups.

The space of admissible equivariant bilinear maps has dimension $q^4$. A single channel has $3q^2$ parameters, and two independent scalings preserve its value. We prove that these scalings account for the generic fiber dimension: the affine closure $X_{s,t}$ of the one-channel image satisfies
\[
\dim X_{s,t}=3q^2-2
\]
for every $q,s,t$ (Theorem~\ref{thm:one-channel-dim}). The proof uses a permutation Jacobian minor at $(I,I,e_o)$ and is independent of the input types. Consequently, the generic fusion width is at least $\lceil q^4/(3q^2-2)\rceil$.

We next prove that
\[
\dim\overline{\im\Phi_2}=6q^2-4\qquad(q\ge3)
\]
for every central-character pair (Theorem~\ref{thm:second-secant}). We construct two parameter points whose differential images have dimension $3q^2-2$ and intersect trivially. The argument uses the supports of matrix entries and explicit character factors. A separate final step handles characteristic two.

Theorem~\ref{thm:subgroup-lift} is stated for finite abelian groups and arbitrary complex-valued bicharacters. If an $r$-channel parameter point on a subgroup $H\le K$ has differential rank $r(3|H|-2)$ and every parameter vector is nonzero, then the differential at its zero extension to $K$ has rank $r(3|K|-2)$. In particular, the affine closure of the one-channel image has dimension $3|K|-2$ for any such group and bicharacters (Corollary~\ref{cor:abelian-one-channel}). Exact modular certificates determine all secant dimensions for $q=2,3,4,5$, with generic widths $2,4,6,9$ (Theorem~\ref{thm:sharp-smallq}). Lifting parameter points below the filling order proves the expected dimensions through orders $3,5,8$ for $q=3^f,4^f,5^f$, respectively, for every $f\ge1$ and every type pair (Corollary~\ref{cor:extension-families}). The characteristic-two conclusion also holds for $q=2^f$ with every $f\ge4$, using an additive subgroup of $\F_q\times\{0\}$. These ranges do not reach the first filling order for large extension degree.

Proposition~\ref{prop:triality} shows that permuting the three central weights $(s,t,-s-t)$ preserves channel decompositions up to a linear change of coordinates. This reduces the number of type pairs in the exact checks. Two rank certificates per type class suffice to determine all secant dimensions in each of these small fields by Lemma~\ref{lem:prefix-compression}. For arbitrary $q$, rank factorization also gives a constructive upper bound of $q^2$ channels for every target.

The untwisted case is related to group-based phylogenetic models. In Fourier coordinates its one-channel variety is the toric tripod for $G/Z\cong\F_q^2$. Its dimension $3q^2-2$ already follows from the claw-dimension computation in the proof of Proposition~5.2 of Ba{\~n}os et al.~\citep{BanosEtAl2019}. Theorem~\ref{thm:one-channel-dim} gives a single permutation minor for all input types, including the twisted cases. Untwisted all-orders nondefectivity is a special case of their Conjecture~5.1, and their Proposition~5.3 proves the untwisted $q=2$ second secant. Their Section~5.2 also uses finite-field Jacobian ranks to certify dimensions. Our contributions concern the uniform second-secant theorem, subgroup lifting, and the particular small-field certificates and their consequences.

There are also connections with structured matrix decompositions and equivariant bilinear complexity. Couveignes and Ezome~\citep{CouveignesEzome2023} study the complexity of multiplication in finite field extensions under equivariance. Huhtanen et al.~\citep{HuhtanenEtAl2026} study sums of diagonal--circulant--diagonal matrices. Their Theorem~2.1 gives the rank-one inclusion, and their Definition~4.3 introduces DCD-rank. After replacing the second Fourier index by its negative, our untwisted parameterization has the same form for the additive group $\F_q^2$. This group is not cyclic, so our results do not settle the cyclic DCD-rank question. The rank-factorization upper bound is a standard consequence of rank-one inclusion. Lin et al.~\citep{LinEtAl2025} approximate $SO(3)$-equivariant tensor products by low-rank decompositions, obtaining approximately equivariant networks with a bound on the equivariance error. Our width concerns exact representation by the specified fusion channels.

In \emph{Algebraic Expressivity Certificates for Shallow Polynomial Neural Networks}, Akbari and Jamshidi (2026; \href{https://arxiv.org/abs/2609.28500}{arXiv:2609.28500}, Propositions~2--3 and Theorem~6) study Veronese secant closures, elimination certificates and the quadratic determinantal case for bias-free shallow power-sum networks. Their width counts scalar pure powers; our width counts the prescribed equivariant bilinear channels. Remark~\ref{rem:waring-comparison} gives a one-channel example whose scalar basepoint polynomial has quadratic Waring rank $2q^2$. The distinction between an image and its closure is relevant to both parameterizations.

Sections 2--4 describe the induced representations, their commutants and the admissible bilinear maps. Sections 5--9 contain the dimension and width results, the type symmetries and the exact certificates. Section~\ref{sec:secant-geometry} gives the toric interpretation and distinguishes fusion width from quadratic Waring rank. Sections~\ref{sec:width-conjecture}--\ref{sec:conclusion} discuss the remaining conjecture, reproducibility and further questions. Appendices~\ref{app:polynomial} and~\ref{app:approximation} give polynomial realization and approximation consequences and explain their relation to earlier results.

\section{Finite Heisenberg groups and induced representations}

All representation spaces and algebraic varieties are over $\C$. Finite coefficient fields are used only for computational verification.

Let \(q=p^f\) be a prime power and write
\[
g(a,b,c)=
\begin{pmatrix}
1&a&b\\0&1&c\\0&0&1
\end{pmatrix}.
\]
Then
\begin{equation}
 g(a,b,c)g(a',b',c')
 =g(a+a',\,b+b'+ac',\,c+c').
 \label{eq:group-law}
\end{equation}
The center and quotient are
\[
Z=\{g(0,b,0):b\in\F_q\}\cong(\F_q,+),
\qquad
G/Z\cong(\F_q^2,+).
\]
Fix a nontrivial additive character \(\psi:\F_q\to\C^\times\).  The nondegenerate trace pairing implies that
\[
\psi_t(b):=\psi(tb),\qquad t\in\F_q,
\]
exhausts all additive characters of \(\F_q\).

\begin{definition}[Induced representation]
For \(t\in\F_q\), define
\[
\mathcal H_t=\Ind_Z^G\psi_t
=
\{F:G\to\C:F(gz)=\psi_t(z)^{-1}F(g)\}.
\]
The left action is \([L_hF](g)=F(h^{-1}g)\).
\end{definition}

Choose the section
\[
r(x,y)=g(x,0,y),\qquad (x,y)\in\F_q^2.
\]
Restriction to the section identifies \(\mathcal H_t\) with \(\C^{\F_q^2}\), hence
\begin{equation}
\dim_{\C}\mathcal H_t=q^2.
\label{eq:Hdim}
\end{equation}
A direct calculation from \eqref{eq:group-law} gives, for \(h=g(\alpha,\beta,\gamma)\),
\begin{equation}
[L_h^tf](x,y)
=
\psi_t(\beta-\alpha\gamma+\alpha y)
 f(x-\alpha,y-\gamma).
\label{eq:left-action}
\end{equation}
In particular, the center acts on all of \(\mathcal H_t\) by the scalar character \(\psi_t\).

For later use we record the irreducible representations explicitly.

\begin{proposition}[Elementary finite Stone--von Neumann decomposition]
For each \(t\neq0\), define \(\pi_t\) on \(\C^{\F_q}\) by
\begin{equation}
[\pi_t(g(a,b,c))\varphi](x)
=
\psi_t(b+cx)\varphi(x+a).
\label{eq:schrodinger}
\end{equation}
Then \(\pi_t\) is irreducible of dimension \(q\), has central character \(\psi_t\), and the representations \(\pi_t\), \(t\neq0\), are pairwise inequivalent.  Together with the \(q^2\) one-dimensional characters factoring through \(G/Z\), they exhaust the irreducible representations of \(G\).
\end{proposition}

\begin{proof}
Equation \eqref{eq:group-law} shows directly that \eqref{eq:schrodinger} is a representation.  Its character is
\[
\chi_t(g(a,b,c))=
\begin{cases}
q\psi_t(b),&a=c=0,\\
0,&\text{otherwise}.
\end{cases}
\]
Indeed, a nonzero translation has zero trace, while for \(a=0\), \(c\neq0\), the trace is
\(\psi_t(b)\sum_x\psi_t(cx)=0\).
Therefore
\[
\langle\chi_t,\chi_t\rangle_G
=\frac1{q^3}\sum_{b\in\F_q}q^2=1,
\]
so \(\pi_t\) is irreducible.  Distinct \(t\) give distinct central characters and hence inequivalent representations.  Since \([G,G]=Z\), there are \(q^2\) one-dimensional representations.  Their squared dimensions contribute \(q^2\); the \(q-1\) representations \(\pi_t\) contribute \((q-1)q^2\).  The sum is \(q^3=|G|\), so the list is complete.
\end{proof}

\begin{corollary}
For \(t\neq0\),
\[
\mathcal H_t\cong q\,\pi_t,
\]
while \(\mathcal H_0\) is the regular representation of \(G/Z\).
\end{corollary}

\begin{proof}
For \(t\neq0\), Frobenius reciprocity gives
\[
[\mathcal H_t:\pi_t]
=
\dim\Hom_Z(\pi_t|_Z,\psi_t)=q,
\]
because \(Z\) acts on the \(q\)-dimensional space of \(\pi_t\) by \(\psi_t\).  No other irreducible has this central character.  The statement for \(t=0\) follows because \(\Ind_Z^G\mathbf1\) is the regular representation of the quotient.
\end{proof}

\section{Twisted convolutions and commutants}

For \((u,v)\in\F_q^2\), define the twisted shift
\begin{equation}
[R_{u,v}^tf](x,y)
=
\psi_t(-xv)f(x+u,y+v).
\label{eq:twisted-shift}
\end{equation}
These are right translations expressed in the chosen section.

\begin{lemma}[Commutation and cocycle law]
For all \(h\in G\) and \((u,v)\in\F_q^2\),
\[
L_h^tR_{u,v}^t=R_{u,v}^tL_h^t.
\]
Moreover,
\begin{equation}
R_{u,v}^tR_{u',v'}^t
=
\psi_t(-uv')R_{u+u',v+v'}^t.
\label{eq:cocycle}
\end{equation}
\end{lemma}

\begin{proof}
Substitution into \eqref{eq:left-action} and \eqref{eq:twisted-shift} proves commutation.  For the product,
\[
\psi_t(-xv)\psi_t(-(x+u)v')
=
\psi_t(-uv')\psi_t(-x(v+v')),
\]
which yields \eqref{eq:cocycle}.
\end{proof}

\begin{theorem}[Commutants and twisted convolutions]
For every \(t\in\F_q\),
\begin{equation}
\End_G(\mathcal H_t)
=
\Span_{\C}\{R_{u,v}^t:(u,v)\in\F_q^2\},
\qquad
\dim_{\C}\End_G(\mathcal H_t)=q^2.
\label{eq:commutant}
\end{equation}
Consequently every \(G\)-equivariant linear map \(T:\mathcal H_t\to\mathcal H_t\) has a unique twisted-convolution representation
\begin{equation}
[T_kf](x,y)
=
\sum_{u,v\in\F_q}
 k(u,v)\psi_t(-xv)f(x+u,y+v).
\label{eq:twisted-convolution}
\end{equation}
\end{theorem}

\begin{proof}
If \(t=0\), then \(\mathcal H_0\) is the regular representation of the abelian group \(G/Z\), hence the multiplicity-free sum of its \(q^2\) characters.  Therefore \(\dim\End_G(\mathcal H_0)=q^2\).  If \(t\neq0\), then \(\mathcal H_t\cong q\pi_t\), so Schur's lemma gives
\[
\End_G(\mathcal H_t)\cong M_q(\C),
\qquad
\dim\End_G(\mathcal H_t)=q^2.
\]
The twisted shifts lie in the commutant by the previous lemma.  They are linearly independent: applying a linear relation to the delta function at \((0,0)\) produces functions supported at distinct points \((-u,-v)\).  Since there are \(q^2\) shifts, they form a basis.
\end{proof}

\begin{remark}
For \(t=0\), \eqref{eq:twisted-convolution} is ordinary convolution on the abelian quotient \(\F_q^2\).  For \(t\neq0\), the phase in \eqref{eq:twisted-convolution} records the Heisenberg cocycle.  Thus the same coordinate space \(\F_q^2\) carries different convolution algebras according to the central character.
\end{remark}

\section{Equivariant bilinear maps}

Throughout, ``bilinear'' means complex-bilinear.  For types \(s,t\), pointwise multiplication
\[
m_{s,t}(f,g)=fg
\]
is a \(G\)-equivariant bilinear map
\(
\mathcal H_s\times\mathcal H_t\to\mathcal H_{s+t}
\), because the phases in \eqref{eq:left-action} multiply.

\begin{theorem}[Selection rule and exact dimension]
For \(s,t,u\in\F_q\),
\begin{equation}
\Bil_G(\mathcal H_s,\mathcal H_t;\mathcal H_u)=0
\qquad\text{if }u\neq s+t.
\label{eq:selection-zero}
\end{equation}
If \(u=s+t\), then
\begin{equation}
\dim_{\C}\Bil_G(\mathcal H_s,\mathcal H_t;\mathcal H_{s+t})=q^4.
\label{eq:bil-dim}
\end{equation}
\end{theorem}

\begin{proof}
Complex-bilinear equivariant maps are naturally
\[
\Hom_G(\mathcal H_s\otimes\mathcal H_t,\mathcal H_u).
\]
For finite groups induction agrees with coinduction. Since \(\mathcal H_u=\Ind_Z^G\psi_u\), Frobenius reciprocity in the target variable gives
\begin{equation}
\Hom_G(\mathcal H_s\otimes\mathcal H_t,\mathcal H_u)
\cong
\Hom_Z((\mathcal H_s\otimes\mathcal H_t)|_Z,\psi_u).
\label{eq:frob-bil}
\end{equation}
The center acts on the tensor product by \(\psi_s\psi_t=\psi_{s+t}\).  Since
\(\dim(\mathcal H_s\otimes\mathcal H_t)=q^4\), as a \(Z\)-module it is
\(\psi_{s+t}^{\oplus q^4}\).  Equation \eqref{eq:frob-bil} therefore has dimension zero unless \(u=s+t\), and dimension \(q^4\) in the admissible case.
\end{proof}

Let \(N=q^2\) and fix the basepoint \(o=(0,0)\in G/Z\).

\begin{lemma}[Basepoint functional model]\label{lem:basepoint-functional}
For every admissible pair \((s,t)\), evaluation at \(o\) induces a linear isomorphism
\[
\operatorname{ev}_o:\Bil_G(\mathcal H_s,\mathcal H_t;\mathcal H_{s+t})
\longrightarrow \Bil_{\C}(\mathcal H_s,\mathcal H_t;\C).
\]
Equivalently, every scalar complex-bilinear functional at the basepoint extends uniquely to a \(G\)-equivariant \(\mathcal H_{s+t}\)-valued bilinear map.
\end{lemma}

\begin{proof}
Apply Frobenius reciprocity to \(W=\mathcal H_s\otimes\mathcal H_t\):
\[
\Hom_G(W,\Ind_Z^G\psi_{s+t})
\cong
\Hom_Z(W,\psi_{s+t}).
\]
Under the induced-function realization, this isomorphism is evaluation of the output at the identity coset \(o=Z\).  On \(W\), the center acts by the scalar character
\(\psi_s\psi_t=\psi_{s+t}\).  Consequently every complex-linear functional \(W\to\C\) is automatically \(Z\)-equivariant when the target \(\C\) is given the character \(\psi_{s+t}\).  Hence
\[
\Hom_Z(W,\psi_{s+t})=\Hom_{\C}(W,\C)
\cong \Bil_{\C}(\mathcal H_s,\mathcal H_t;\C),
\]
which proves both existence and uniqueness of the equivariant extension.
\end{proof}

Choosing the standard coordinate basis of \(\C^{\F_q^2}\), Lemma~\ref{lem:basepoint-functional} yields
\begin{equation}
\Bil_G(\mathcal H_s,\mathcal H_t;\mathcal H_{s+t})
\cong M_N(\C),
\label{eq:basepoint-matrix}
\end{equation}
where \(M\in M_N(\C)\) represents the bilinear functional
\(B(f,g)(o)=f^TMg\).  We use this identification in the Jacobian computations below.

\section{Fusion width and a general lower bound}\label{sec:fusion-lower}

\begin{definition}[Fusion channel]
A fusion channel from \(\mathcal H_s\times\mathcal H_t\) to \(\mathcal H_{s+t}\) is
\begin{equation}
\Gamma_{A,B,C}(f,g)
=
C[(Af)(Bg)],
\label{eq:channel}
\end{equation}
where
\[
A\in\End_G(\mathcal H_s),\quad
B\in\End_G(\mathcal H_t),\quad
C\in\End_G(\mathcal H_{s+t}).
\]
For \(r\ge1\), define the polynomial parameterization
\begin{equation}
\Phi_r((A_\ell,B_\ell,C_\ell)_{\ell=1}^r)
=
\sum_{\ell=1}^r\Gamma_{A_\ell,B_\ell,C_\ell}.
\label{eq:Phi}
\end{equation}
\end{definition}

The operators $A$ and $B$ act on the two inputs, and $C$ acts on their pointwise product. If $d$ is the basepoint row of $C$, the coefficient matrix of one channel is
\[
M=A^T\operatorname{diag}(d)B,
\qquad A\in\End_G(\mathcal H_s),\quad B\in\End_G(\mathcal H_t),\quad d\in\C^N.
\]
The vector $d$ is arbitrary because the basepoint rows of the twisted shifts form the coordinate basis. Thus fusion width is the number of summands in this structured matrix decomposition. It is not ordinary matrix rank: taking $A=B=I$ and $d=\mathbf 1$ gives $M=I_N$ in one channel. The same example is considered in Remark~\ref{rem:waring-comparison}.

\begin{definition}[Generic fusion width]
Let
\[
V_{s,t}=\Bil_G(\mathcal H_s,\mathcal H_t;\mathcal H_{s+t}).
\]
The width of an individual target is the least nonnegative number of channels representing it. We set \(\operatorname{width}(0)=0\) and \(\Phi_0=0\). The generic fusion width \(r_{\rm gen}(q;s,t)\) is the least \(r\) for which \(\Phi_r\) is dominant onto \(V_{s,t}\).
\end{definition}

Here and throughout, \emph{generic} is meant in the complex-algebraic sense: the asserted property holds on a nonempty Zariski-open subset of the affine target space.  Consequently its complement lies in a proper algebraic subset and has Lebesgue measure zero after identifying the target with a real Euclidean space.  Because the image of a morphism is constructible, dominance here has an exact representability meaning, not merely a border-rank meaning.

\begin{proposition}[Dominance implies exact generic width]\label{prop:dominance-exact}
If \(\Phi_r\) is dominant and \(\Phi_{r-1}\) is not, then there exists a nonempty Zariski-open subset \(U\subset V_{s,t}\) such that every element of \(U\) has exact fusion width \(r\).
\end{proposition}

\begin{proof}
By Chevalley's theorem, the dense constructible set \(\im\Phi_r\) contains a nonempty open subset of the irreducible affine space \(V_{s,t}\).  Since \(\overline{\im\Phi_{r-1}}\) is a proper closed subset, removing it leaves a nonempty open subset on which width is at least \(r\) and at most \(r\).
\end{proof}

\begin{proposition}[Constructive exact upper bound]\label{prop:constructive-upper}
Let $M$ be the basepoint coefficient matrix of an admissible equivariant bilinear map $B$. Then
\[
\operatorname{width}(B)\le \operatorname{rank}(M)\le N=q^2.
\]
Given a rank factorization of $M$, the corresponding channel kernels are obtained explicitly. A factorization and the resulting kernels can be computed using $O(N^3)$ arithmetic operations.
\end{proposition}
\begin{proof}
Write $M=\sum_{\ell=1}^{k}a_\ell b_\ell^T$, where $k=\operatorname{rank}(M)$. The basepoint row of the twisted convolution with coefficient kernel $a_\ell$ is $a_\ell^T$, independently of its type; similarly for $b_\ell$. Choose the output kernel $c=e_o$, so that its convolution operator is $C=I$ and its basepoint row is $e_o^T$. Consequently
\[
 A(a_\ell)^T\operatorname{diag}(e_o)B(b_\ell)=a_\ell b_\ell^T.
\]
Each summand is therefore one channel. Evaluation at the basepoint is an isomorphism on admissible equivariant bilinear maps, so equality of these matrices gives equality of the maps. Gaussian rank factorization costs $O(N^3)$ arithmetic operations, and writing the kernels costs $O(kN)$ operations. This is an arithmetic, rather than bit, complexity statement.
\end{proof}

\begin{corollary}[A worst-case consequence of dominance]
If $\Phi_r$ is dominant, then every target has fusion width at most $\min\{N,2r\}$. In particular, the worst-case widths for $q=2,3,4,5$ are at most $4,8,12,18$, respectively.
\end{corollary}
\begin{proof}
The image contains a nonempty Zariski-open subset $U$ of $V_{s,t}$. For any $M\in V_{s,t}$, the two nonempty open subsets $U$ and $M-U$ intersect because the affine space $V_{s,t}$ is irreducible. Thus $M$ is a sum of two elements in $\operatorname{im}\Phi_r$. Apply Proposition~\ref{prop:constructive-upper} for the other bound.
\end{proof}

Proposition~\ref{prop:constructive-upper} gives a decomposition of every target, but need not attain the generic minimum. The factor-two bound follows from the standard dense-open-set argument; see Blekherman and Teitler~\citep{BlekhermanTeitler2015}. The rank-one construction also occurs in the DCD-matrix setting of Huhtanen et al.~\citep{HuhtanenEtAl2026}.

\begin{lemma}[Scaling symmetries]\label{lem:scaling-bound}
The one-channel parameterization is invariant under
\begin{equation}
(A,B,C)\longmapsto
(\lambda A,\mu B,(\lambda\mu)^{-1}C),
\qquad
\lambda,\mu\in\C^\times.
\label{eq:gauge}
\end{equation}
Hence the Zariski closure of the one-channel image has dimension at most
\[
3q^2-2.
\]
More generally,
\begin{equation}
\dim\overline{\im\Phi_r}\le r(3q^2-2).
\label{eq:dim-upper}
\end{equation}
\end{lemma}

\begin{proof}
Each commutant has dimension \(q^2\).  On the dense open subset where \(A,B,C\neq0\), the torus \((\C^\times)^2\) in \eqref{eq:gauge} acts freely and preserves every channel.  For \(r\) channels the product torus \((\C^\times)^{2r}\) acts freely on a dense open subset of the parameter space and is contained in each fiber.  The fiber-dimension theorem gives \eqref{eq:dim-upper}.
\end{proof}

We now prove that the preceding dimension bound is attained for every prime power and every type pair.

\begin{theorem}[Dimension of the one-channel variety]\label{thm:one-channel-dim}
Let $X_{s,t}\subset V_{s,t}$ be the affine Zariski closure of the one-channel image. For every prime power $q$ and every $s,t\in\F_q$,
\[
\dim X_{s,t}=3q^2-2.
\]
The projectivized variety has dimension $3q^2-3$, and the generic fiber has dimension two. At $(A,B,d)=(I,I,e_o)$, the Jacobian has an explicitly selected $(3q^2-2)$-square permutation submatrix. Its determinant is $(-1)^{q^2-1}$ in the ordering below, and all its singular values equal one.
\end{theorem}
\begin{proof}
Put $N=q^2$ and $K=\F_q^2$. In the basepoint matrix model, the one-channel map is
\[
F(A,B,d)=A^T\operatorname{diag}(d)B.
\]
The vector $d$ is arbitrary, since the basepoint rows of the twisted shifts are the coordinate rows. At $(I,I,d)$ the differential is
\begin{equation}
 dF_{(I,I,d)}(X,Y,z)=X^T\operatorname{diag}(d)
 +\operatorname{diag}(z)+\operatorname{diag}(d)Y.
 \label{eq:one-channel-differential}
\end{equation}
Write $X=\sum_h x_hR_h^s$ and $Y=\sum_h y_hR_h^t$. For $h=(u,v)\ne o$ and $i=(i_x,i_y)$, the $(i,i-h)$ entry is
\begin{equation}
 x_h\psi_s(-(i_x-u)v)d_{i-h}
 +y_{-h}\psi_t(i_xv)d_i.
 \label{eq:offdiag-kernel}
\end{equation}
Now set $d=e_o$. Select the $N$ columns $z_i$ and the $N$ diagonal rows. For each $h\ne o$, select columns $x_h,y_{-h}$ and rows $(o,-h),(h,o)$, in that order. All selected rows and columns are distinct, including in characteristic two. Since the basepoint row of $R_h$ is $e_h^T$, the selected matrix is
\begin{equation}
 I_N\oplus\bigoplus_{h\ne o}
 \begin{pmatrix}0&1\\1&0\end{pmatrix}.
 \label{eq:explicit-block-minor}
\end{equation}
This proves the determinant and singular-value assertions, and gives differential rank at least $3N-2$. The dimension bound in Lemma~\ref{lem:scaling-bound} gives the reverse inequality. Since $X_{s,t}$ is a nonzero affine cone, its projectivization has dimension $3N-3$. The generic-fiber dimension is $3N-(3N-2)=2$ by the fiber-dimension theorem.
\end{proof}

\begin{remark}[Scope of the explicit minor]
The selected minor has spectral condition number one in the stated coordinates. This condition number concerns the differential at the specified point. A generic fiber of dimension two need not be a single scaling orbit: common twisted shifts also give discrete parameter symmetries.
\end{remark}

\begin{theorem}[Lower bound for generic width]\label{thm:general-lower}
For every prime power \(q\) and every \(s,t\in\F_q\),
\begin{equation}
 r_{\rm gen}(q;s,t)
 \ge
 \left\lceil\frac{q^4}{3q^2-2}\right\rceil.
\label{eq:general-lower}
\end{equation}
\end{theorem}

\begin{proof}
The target dimension is \(q^4\) by \eqref{eq:bil-dim}.  By \eqref{eq:dim-upper}, the closure of the $r$-channel image has dimension at most $r(3q^2-2)$.  If \(r(3q^2-2)<q^4\), dominance is impossible.  Solving for \(r\) proves the claim.
\end{proof}

\begin{remark}
By Theorems~\ref{thm:one-channel-dim} and~\ref{thm:second-secant}, any failure of the expected secant dimension must occur at order at least three.
\end{remark}

\section{Nondefectivity of the second secant}\label{sec:second-secant}

\begin{theorem}[Dimension of the second secant]\label{thm:second-secant}
Let $q\ge3$ be a prime power and let $s,t\in\F_q$ be arbitrary. Then
\[
 \dim\overline{\operatorname{im}\Phi_2}=6q^2-4.
\]
Thus the second projective secant of the one-channel variety has its expected
dimension $6q^2-5$. Together with the $q=2$ case of Theorem~\ref{thm:sharp-smallq}, this gives
\[
 \dim\overline{\operatorname{im}\Phi_2}
 =\min\{6q^2-4,q^4\}
\]
for every prime power $q$ and every pair of central types. The argument below applies to every $q\ge3$; only the case $q=2$ uses a rank certificate.
\end{theorem}
\begin{proof}
Write $K=\F_q^2$, $N=q^2$, $o=(0,0)$, $a=(1,0)$ and $b=(0,1)$. Put
\[
 H=\{o,a,-b,a-b\},\qquad
 J=\begin{cases}K\setminus\{o\},&\operatorname{char}\F_q\ne2,\\
 K\setminus H,&\operatorname{char}\F_q=2.
 \end{cases}
\]
The set $H$ has four elements. It is a subgroup only in the characteristic-two
case. A stripe of difference $h$ of a matrix consists of its entries
$M_{i,i-h}$, $i\in K$. Distinct stripes have disjoint supports.

Choose $d\in\C^K$ satisfying these two conditions:
\begin{enumerate}
\item[(i)] For every $h=(u,v)\notin H$, the following two vectors, restricted
to $K\setminus E_h$ with $E_h=\{o,a,h,h+b\}$, are linearly independent:
\[
 \left(\psi_s(-(i_x-u)v)d_{i-h}\right)_i,
 \qquad \left(\psi_t(i_xv)d_i\right)_i.
\]
\item[(ii)] The three vectors on $J$ given by
\[
 D=(d_k)_k,\quad D_a=(d_{k+a})_k,\quad
 W=(\psi_{s+t}(k_x)d_{k+b})_k
\]
are linearly independent.
\end{enumerate}
We first prove that such a vector $d$ exists. Since
$N\ge9$, the complement of $E_h$ has at least $N-4\ge5$ elements.
For any two distinct indices $i,j$ in that complement, the determinant of
the corresponding two rows in (i) is a difference of nonzero scalar multiples
of the monomials $d_{i-h}d_j$ and $d_i d_{j-h}$. These monomials are distinct:
equality of the multisets of their indices would force either $h=o$ or $i=j$.
Hence the selected determinant is a nonzero polynomial.
For (ii) in odd characteristic, choose $w\notin\{o,a,b\}$; in characteristic
two choose $w\notin H$, which is possible since $q\ge4$ and $N\ge16$.
Specializing $d=e_w$, the three distinct indices $w,w-a,w-b$ all belong
to $J$. The columns $D,D_a,W$ on these rows form a diagonal matrix with
nonzero diagonal (up to their ordering). Thus this three-square determinant
is also a nonzero polynomial. The simultaneous nonvanishing conditions are
finitely many nonempty principal Zariski-open subsets of the irreducible
space $\C^K$, and their intersection is nonempty.

Choose complex scalars
\[
 (\alpha,\beta)=\begin{cases}(1,1),&\operatorname{char}\F_q\ne2,\\
 (2,3),&\operatorname{char}\F_q=2,
 \end{cases}
\]
and consider the two parameter points
\[
 P=(I,I,d),\qquad Q=(A,B,e_o),\qquad
 A=I+\alpha R_a^s,\quad B=I+\beta R_b^t.
\]
Both $A$ and $B$ are invertible. The shifts $R_a^s,R_b^t$ have order $p$,
where $p=\operatorname{char}\F_q$. For odd $p$, $-1$ is not an eigenvalue.
For $p=2$ their squares are the identity and $\alpha^2,\beta^2\ne1$.
Let $T_P,T_Q$ denote the two differential images. At $P$ the diagonal stripe has dimension $N$, and every stripe outside $H$ has dimension two by (i). For the stripes $a,-b,a-b$, evaluate formula~\eqref{eq:offdiag-kernel} at $(k+a,k)$, $(k,k+b)$ and $(k+a,k+b)$, respectively, with $k\in J$. This gives the pairs
\[
 (D,D_a),\qquad (W,D),\qquad (W,\psi_t(-1)D_a),
\]
after multiplication of rows by nonzero scalars. In the last two pairs the row indexed by $k$ is divided by $\psi_t(-k_x)$. Condition (ii) makes each pair independent. Since the stripes have disjoint supports,
\[
 \dim T_P=N+2(N-1)=3N-2.
\]
Let $T_o$ denote the differential image at $(I,I,e_o)$; it consists exactly
of the matrices supported on the diagonal together with row $o$ and column
$o$. The changes $X\mapsto XA$, $Y\mapsto YB$ are bijections of the
respective commutants. Differentiation therefore gives
\[
 T_Q=A^TT_oB,\qquad \dim T_Q=3N-2.
\]

We show that $T_P\cap T_Q=\{0\}$. For an element $M$ in this intersection,
write
\[
 M=xv^T+uy^T+A^T\operatorname{diag}(z)B,\qquad
 u=e_o+\alpha e_a,\quad v=e_o+\beta e_b,
\]
where $x,y,z\in\C^K$. The vectors $x,y$ are arbitrary because the basepoint
row maps of the two commutants are isomorphisms. The term indexed by $k$
in $A^T\operatorname{diag}(z)B$ has positions
\[
 (k,k),\quad(k+a,k),\quad(k,k+b),\quad(k+a,k+b),
\]
so this term is supported only on stripes $H$. For $h\notin H$, the first
two terms may be nonzero only at rows $i\in E_h$. On the other hand,
$M\in T_P$ restricts on that stripe to a linear combination of the two
vectors in (i). Vanishing on $K\setminus E_h$ forces its two coefficients
to be zero. Thus all stripes outside $H$ vanish.

For $i\notin\{o,a\}$, the entries $(i,o)$ and $(i,b)$ of the first two
terms are $x_i$ and $\beta x_i$. Unless both differences $i,i-b$ belong
to $H$, a vanishing stripe forces $x_i=0$. Similarly, entries $(o,j),(a,j)$
for $j\notin\{o,b\}$ are $y_j,\alpha y_j$. Consequently,
\[
 \operatorname{supp}x\subseteq H\cap(H+b),\qquad
 \operatorname{supp}y\subseteq(-H)\cap(a-H).
\]
In odd characteristic these two intersections are $\{o,a\}$ and
$\{o,b\}$, respectively. In characteristic two both intersections equal
$H$. These identities follow by comparing the independent coordinates
along $a,b$; they also hold in the smallest odd characteristic $p=3$.

For $k\in J$, the first two terms therefore contribute nothing at the
three positions $(k+a,k),(k,k+b),(k+a,k+b)$. Formula~\eqref{eq:offdiag-kernel} for $T_P$
and division by the nonzero phase $\psi_t(-k_x)$ in the last two positions
give constants $\lambda_j,\mu_j$ such that
\begin{align*}
 \alpha z_k&=\lambda_1d_k+\mu_1d_{k+a},\\
 \beta z_k&=\lambda_2\psi_{s+t}(k_x)d_{k+b}+\mu_2d_k,\\
 \alpha\beta z_k&=\lambda_3\psi_{s+t}(k_x)d_{k+b}
                 +\mu_3\psi_t(-1)d_{k+a}.
\end{align*}
Thus $(z_k)_{k\in J}$ lies in the intersection of the three planes
\[
 \operatorname{span}(D,D_a),\quad
 \operatorname{span}(W,D),\quad
 \operatorname{span}(W,D_a).
\]
Their intersection is zero by (ii). Hence $z_k=0$ for $k\in J$ and each
pair of coefficients $\lambda_j,\mu_j$ is zero. These are the coefficients
of the full stripes of $M\in T_P$, so the stripes $a,-b,a-b$ vanish
at every index, not just on $J$. We conclude that all off-diagonal
entries of $M$ are zero.

It remains to treat the possible diagonal entries. In odd characteristic,
only $z_o$ remains. The matrix $A^T\operatorname{diag}(z)B=z_o uv^T$ can
be absorbed into $xv^T+uy^T$. This matrix is supported on
$\{o,a\}\times\{o,b\}$; its only possible diagonal entry is $M_{oo}$.
Because $\alpha=\beta=1$, its entries obey
\[
 M_{oo}+M_{ab}-M_{ob}-M_{ao}=0.
\]
The other three entries are off-diagonal and vanish. Hence $M_{oo}=0$.

In characteristic two, both $xv^T+uy^T$ and the remaining terms indexed by
$k\in H$ have support in $H\times H$. Thus
$M=\operatorname{diag}(m)$ with $\operatorname{supp}m\subseteq H$.
Since $T_Q=A^TT_oB$ and both shifts square to the identity,
\[
 L=(I-\alpha(R_a^s)^T)M(I-\beta R_b^t)\in T_o;
\]
the scalar denominators omitted from $A^{-T}$ and $B^{-1}$ are nonzero.
The four terms in this expression have the four distinct differences
$o,a,b,a+b$. Put $c=a+b$ and $\theta_k=\psi_t(-k_x)$. Direct evaluation
at four off-diagonal positions gives
\[
 L_{a,b}=\alpha\beta m_o,\quad
 L_{b,a}=\alpha\beta\theta_c m_c,\quad
 L_{c,b}=-\alpha m_b,\quad
 L_{a,c}=-\beta\theta_a m_a.
\]
These entries avoid row $o$, column $o$, and the diagonal, so $L\in T_o$
forces each to vanish. All four scalar multipliers are nonzero; therefore
$m_o=m_a=m_b=m_c=0$.

In both cases $M=0$, proving the required trivial intersection. The
differential of $\Phi_2$ at $(P,Q)$ has rank $6N-4$, which equals the
universal upper bound. For $N\ge9$ this is less than $N^2$, so the
asserted non-filling second-secant dimension follows. For $q=2$ the exact certificate in Section~\ref{sec:certificates} supplies the filling dimension $16$.
\end{proof}

\begin{remark}
The admissible vectors $d$ form a nonempty Zariski-open set defined by the two-square and three-square minors used in the proof. The higher-secant question is considered in Section~\ref{sec:subgroup-lift} and in Conjecture~\ref{conj:general-fusion}.
\end{remark}

\section{Symmetries of central-character pairs}\label{sec:type-symmetry}

In this section, we describe the changes of type that preserve fusion width and secant dimensions.

\begin{proposition}[Scaling and exchange symmetries]
Fusion width is invariant under
\[
(s,t)\longmapsto(\lambda s,\lambda t),\qquad \lambda\in\F_q^\times,
\]
and under exchange \((s,t)\leftrightarrow(t,s)\). These transformations induce linear isomorphisms of the ambient bilinear spaces that preserve every channel image and every secant closure.
\end{proposition}

\begin{proof}
For \(\lambda\neq0\), the map
\[
\sigma_\lambda:g(a,b,c)\longmapsto g(\lambda a,\lambda b,c)
\]
is a group automorphism.  On the center it sends \(b\mapsto\lambda b\), and therefore transports central character \(\psi_t\) to \(\psi_{\lambda t}\).  Precomposition with \(\sigma_\lambda\) transports the induced representation spaces, their commutants, pointwise products, and hence every channel decomposition.  Thus simultaneous nonzero scaling preserves fusion width.  Exchange invariance follows from commutativity of pointwise multiplication after interchanging the two preprocessing operators.
\end{proof}

\begin{proposition}[Triality of the three central weights]\label{prop:triality}
The one-channel varieties and all their secant closures are preserved, up to linear isomorphism of their ambient spaces, by permutations of the triple
\[
(s,t,-s-t).
\]
In particular, besides simultaneous scaling and input exchange, the transformation $(s,t)\mapsto(-s-t,t)$ preserves the entire secant-dimension profile and every exact fusion-width stratum.
\end{proposition}
\begin{proof}
For every type $u$, the pairing
\[
\beta_u(f,h)=\sum_{x\in\F_q^2} f(x)h(x),
\qquad f\in\mathcal H_u,\quad h\in\mathcal H_{-u},
\]
is nondegenerate and $G$-invariant: the opposite central phases cancel, and the remaining action permutes the summands. It identifies $\mathcal H_{-u}$ with the dual of $\mathcal H_u$. Consequently matrix transpose induces a linear isomorphism
$\End_G(\mathcal H_u)\to\End_G(\mathcal H_{-u})$.
Pairing an output of type $u=s+t$ with a third input $h$ of type $-u$ identifies each equivariant bilinear map with a $G$-invariant trilinear form. A channel becomes
\[
\beta_u(C[(Af)(Bg)],h)
=\sum_x (Af)(x)(Bg)(x)(C^T h)(x).
\]
This expression is symmetric in the three factors. Permuting the factors and their preprocessing operators, with transpose when a factor is designated as the output, carries channels bijectively to channels. These transformations are linear on the ambient spaces of multilinear forms and preserve sums of any number of channels, proving all assertions.
\end{proof}

For pairs with both entries nonzero, simultaneous scaling reduces to \((1,r)\), and exchange sends \(r\mapsto r^{-1}\).  Table~\ref{tab:types} gives the orbit representatives retained in the exact certificates:
\begin{table}[ht]
\caption{Archived type representatives under scaling and input exchange.}
\label{tab:types}
\centering
\begin{tabular}{c|l}
\toprule
field & representatives of ordered type pairs \\
\midrule
\(\F_2\) & \((0,0),(0,1),(1,1)\)\\
\(\F_3\) & \((0,0),(0,1),(1,1),(1,2)\)\\
\(\F_4\) & \((0,0),(0,1),(1,1),(1,\alpha)\), \(\alpha^2+\alpha+1=0\)\\
\(\F_5\) & \((0,0),(0,1),(1,1),(1,2),(1,4)\)\\
\bottomrule
\end{tabular}

\end{table}
These representatives cover all ordered pairs in the corresponding fields. They retain the original scaling/exchange convention for compatibility with the frozen certificates. Triality gives the smaller lists
\[
\begin{array}{c|l}
q&\text{representatives under scaling and triality}\\\hline
2&(0,0),(0,1)\\
3&(0,0),(0,1),(1,1)\\
4&(0,0),(0,1),(1,\alpha)\\
5&(0,0),(0,1),(1,1).
\end{array}
\]
Indeed, a triple with a zero entry and not all entries zero is equivalent to \((0,1,-1)\). In \(\F_3\) all nonzero zero-sum triples are proportional to \((1,1,1)\); in \(\F_4\) they are permutations of \((1,\alpha,\alpha^2)\) up to scale; in \(\F_5\) they are equivalent to \((1,1,3)\). Thus only \(2,3,3,3\) classes are needed, respectively; the additional archived classes are redundant checks.

\section{Exact computer-assisted certification for \texorpdfstring{$q=2,3,4,5$}{q=2,3,4,5}}\label{sec:certificates}

For each type representative, we give exact Jacobian certificates at the width immediately before filling and at the first filling width. The supplement also checks every smaller width.

\subsection{Jacobian coordinates}

Choose the twisted-shift bases of the three commutants.  Write a channel as kernels \(a,b,c\in\C^N\).  Let \(A(a),B(b),C(c)\) be the corresponding \(N\times N\) matrices.  If \(c_o\) denotes the row of \(C(c)\) producing the output at the basepoint, then under the identification \eqref{eq:basepoint-matrix} the channel coefficient matrix is
\begin{equation}
M(a,b,c)=A(a)^T\operatorname{diag}(c_o)B(b).
\label{eq:channel-matrix}
\end{equation}
Thus \(\Phi_r\) is a polynomial map from \(\A^{3Nr}\) to \(\A^{N^2}\), and its Jacobian has \(N^2\) rows and \(3Nr\) columns.

\subsection{Lifting modular rank to characteristic zero}

\begin{lemma}[Modular maximal-minor certificate]\label{lem:modular-lift}
Let \(R\) be the cyclotomic integer ring containing the character values occurring in a fixed \(q\)-model, and let \(J\) be a Jacobian matrix with entries in \(R\) evaluated at a parameter point in \(R\).  Suppose a ring homomorphism \(R\to k\) to a finite field sends some \(m\times m\) minor of \(J\) to a nonzero element of \(k\).  Then that minor is nonzero over \(\C\), so the complex Jacobian rank is at least \(m\).
\end{lemma}

\begin{proof}
If the determinant were zero in \(R\), every ring homomorphism would send it to zero.  A nonzero reduction therefore proves that the algebraic integer determinant is nonzero.  Embedding the cyclotomic ring into \(\C\) preserves nonzeroness.
\end{proof}

If \(m=N^2\), the differential is surjective.  Since source and target are smooth affine spaces, the parameterization is dominant.  Combined with the lower bound, this proves exact generic width.

\begin{lemma}[Maximal rank passes to subcollections]\label{lem:prefix-compression}
Let $P$ and $V$ be complex affine spaces, let $F:P\to V$ be a polynomial one-channel parameterization over $\C$, let $d=\dim\overline{F(P)}$, and let $\Phi_r$ be the sum of $r$ copies of $F$. Suppose that at a point $(p_1,\ldots,p_m)$,
\[
\operatorname{rank}\,d\Phi_m=md.
\]
Then, for every subset $S\subseteq\{1,\ldots,m\}$, the sum parameterization on the points indexed by $S$ has differential rank $|S|d$. Hence $\dim\overline{\operatorname{im}\Phi_r}=rd$ for $1\le r\le m$.
\end{lemma}
\begin{proof}
Write $T_i=\operatorname{im}dF_{p_i}$. Each $T_i$ has dimension at most $d$, since the rank of the differential is bounded above by the dimension of the image closure in characteristic zero. The hypothesis states
\[
md=\dim(T_1+\cdots+T_m)\le\sum_i\dim T_i\le md.
\]
All inequalities are equalities. Thus each $T_i$ has dimension $d$ and their sum is direct. Every subcollection is therefore a direct sum of dimension $|S|d$. The general sum-dimension upper bound proves equality for the corresponding image closures.
\end{proof}
Applied after modular lifting with \(m=r_*-1\), Lemma~\ref{lem:prefix-compression} shows that a penultimate rank \((r_*-1)(3N-2)\), together with a full-rank certificate at \(r_*\), already proves the entire profile. This reduces the necessary number of rank witnesses to two per type class. For \(q=4,5\), the archive additionally retains and verifies every sub-filling minor as a redundant check; the \(q=2,3\) scripts verify each prefix rank directly.

\subsection{Certified cases}

Table~\ref{tab:certificate-fields} summarizes the exact arithmetic choices.  For \(q=3\) and \(q=5\), the auxiliary prime is chosen so that the finite field contains the needed roots of unity.  For \(q=2,4\), additive characters take values \(\{\pm1\}\), so reduction modulo \(101\) is immediate.
The specialization maps used by the certificates are explicit:
\[
\mathbb Z[\zeta_3]\longrightarrow\F_{109},\qquad \zeta_3\longmapsto45,
\]
where \(45\) has exact order \(3\) modulo \(109\), and
\[
\mathbb Z[\zeta_5]\longrightarrow\F_{101},\qquad \zeta_5\longmapsto36,
\]
where \(36\) has exact order \(5\) modulo \(101\).  For \(q=2\) the nontrivial additive character is \((-1)^z\), and for \(q=4\) the additive character is \((-1)^{\operatorname{Tr}(tz)}\), so all character values lie in \(\mathbb Z\) and we use the ordinary reduction \(\mathbb Z\to\F_{101}\).  The verifier checks the required root orders before constructing the Jacobians.

\begin{table}[ht]
\caption{Exact coefficient fields and ranks at the penultimate and filling widths.}
\label{tab:certificate-fields}
\centering\small\setlength{\tabcolsep}{4pt}
\begin{tabular}{c|c|c|c|c|c}
\toprule
\(q\) & \(N=q^2\) & target \(N^2\) & threshold & coefficient field & ranks \\
\midrule
2 & 4 & 16 & 2 & \(\F_{101}\), \(\omega_2=-1\) & \(10\to16\)\\
3 & 9 & 81 & 4 & \(\F_{109}\), \(\omega_3=45\) & \(75\to81\)\\
4 & 16 & 256 & 6 & \(\F_{101}\), \(\psi_t(z)=(-1)^{\Tr(tz)}\) & \(230\to256\)\\
5 & 25 & 625 & 9 & \(\F_{101}\), \(\omega_5=36\) & \(584\to625\)\\
\bottomrule
\end{tabular}

\end{table}
Here the two ranks displayed are for widths \(r_*-1\) and \(r_*\), with
\[
r_*=
\left\lceil\frac{q^4}{3q^2-2}\right\rceil.
\]
For \(q=2\), a fixed width-two witness gives rank \(10\) at width one and full rank \(16\) at width two for each representative \((0,0),(0,1),(1,1)\).  The corresponding certified \(16\times16\) maximal-minor determinants modulo \(101\) are
\[
87,\qquad61,\qquad37,
\]
respectively.  Thus the smallest field already attains the parameter-count threshold \(\lceil16/10\rceil=2\).

For \(q=3\), the same fixed parameter point gives the complete rank sequence
\[
25,50,75,81
\]
for widths \(1,2,3,4\), respectively, in every type-pair orbit.  The certified \(81\times81\) minor determinants modulo \(109\) for the representatives \((0,0),(0,1),(1,1),(1,2)\) are, respectively,
\[
16,\quad21,\quad71,\quad34.
\]
For \(q=4,5\), the supplement stores the parameter arrays as fixed JSON data, records the pivot columns of an explicit maximal minor, and verifies its determinant without using a pseudorandom number generator.  For \(q=4\), the certified \(256\times256\) determinants modulo \(101\) for \((0,0),(0,1),(1,1),(1,\alpha)\) (with the implementation encoding \(\alpha\) by the integer \(2\)) are
\[
91,\quad49,\quad12,\quad47.
\]
For \(q=5\), the certified \(625\times625\) determinants modulo \(101\) for \((0,0),(0,1),(1,1),(1,2),(1,4)\) are
\[
49,\quad78,\quad43,\quad24,\quad85.
\]
A redundant exact recheck at every sub-filling width gives, for every representative type pair,
\[
\begin{array}{c|l}
q&\text{Jacobian ranks at widths }1,2,\ldots,r_*\\
\hline
2&10,16\\
3&25,50,75,81\\
4&46,92,138,184,230,256\\
5&73,146,219,292,365,438,511,584,625.
\end{array}
\]
Before filling, each rank equals $r(3q^2-2)$. At the first filling width it equals $q^4$. All ranks and determinants are computed by exact modular arithmetic.

\begin{theorem}[Secant dimensions for \(q=2,3,4,5\)]\label{thm:sharp-smallq}
For \(q\in\{2,3,4,5\}\), every central-character pair \((s,t)\), and every \(r\ge1\),
\begin{equation}
\dim\overline{\operatorname{im}\Phi_r}
=
\min\{r(3q^2-2),q^4\}.
\label{eq:smallq-full-profile}
\end{equation}
Equivalently, every projective secant through the first filling order has its expected dimension.  In particular,
\begin{equation}
 r_{\rm gen}(2;s,t)=2,\qquad
 r_{\rm gen}(3;s,t)=4,\qquad
 r_{\rm gen}(4;s,t)=6,\qquad
 r_{\rm gen}(5;s,t)=9.
\label{eq:sharp-smallq-values}
\end{equation}
\end{theorem}

\begin{proof}
Theorem~\ref{thm:one-channel-dim} gives the upper bound
\(\dim\overline{\operatorname{im}\Phi_r}\le\min\{r(3q^2-2),q^4\}\).
For each certified type-pair representative, the exact modular Jacobian ranks displayed above attain this bound at every width through \(r_*\).  By Lemma~\ref{lem:modular-lift}, the same ranks hold as lower bounds in characteristic zero, proving equality in \eqref{eq:smallq-full-profile} up to the filling order.  At \(r=r_*\) the map is dominant; hence it remains dominant for every larger \(r\) by setting the extra channels to zero.  Section~\ref{sec:type-symmetry} transports the result from the representatives to all ordered type pairs.  The width values then follow from Proposition~\ref{prop:dominance-exact} and the general lower bound.
\end{proof}

\begin{remark}[Scope of the computation]
The structure of the induced representations, the twisted-convolution characterization, the bilinear selection rule, the target dimension \(q^4\), the exact one-channel dimension \(3q^2-2\), the general lower bound, and the type-pair symmetry reductions are proved analytically.  Theorem~\ref{thm:second-secant} also proves order-two nondefectivity analytically for every $q\ge3$. Computation supplies the complete small-field profiles and hence the base witnesses used in Corollary~\ref{cor:extension-families}; the lifting theorem itself is analytic.
\end{remark}

\section{Lifting witnesses from additive subgroups}\label{sec:subgroup-lift}

\begin{theorem}[Lifting from a subgroup]
\label{thm:subgroup-lift}
Let $K$ be a finite abelian group, let $H\le K$, and write $N=|K|$ and $n=|H|$.
Let $\beta_s,\beta_t:K\times K\to\C^\times$ be bicharacters. Consider the polynomial channel map
\[
 F_K(a,b,d)_{i,j}=\sum_{z\in K} a_{i-z}b_{j-z}d_z
               \beta_s(z,i-z)\beta_t(z,j-z),
 \qquad i,j\in K,
\]
and define $F_H$ by restricting the indices and bicharacters to $H$.
Suppose that, at $r$ parameter triples $p_\ell=(a_\ell,b_\ell,d_\ell)$ on $H$, all three vectors in each triple are nonzero and
\[
 \operatorname{rank}d\!\left(\sum_{\ell=1}^r F_H\right)_{(p_1,\ldots,p_r)}
 =r(3n-2).
\]
Extend each parameter vector by zero from $H$ to $K$. At the resulting parameter point on $K$,
\[
 \operatorname{rank}d\!\left(\sum_{\ell=1}^r F_K\right)
 =r(3N-2).
\]
Consequently the $r$-channel image closure on $K$ has dimension $r(3N-2)$.
\end{theorem}
\begin{proof}
The two scaling directions
\[
 (a_\ell,0,-d_\ell),\qquad (0,b_\ell,-d_\ell)
\]
lie in the kernel of the restricted Jacobian. They are independent, including across channels. The rank hypothesis says that these $2r$ vectors span its entire kernel.
It follows that each of the three partial Jacobians, obtained by varying only all $a$-vectors, only all $b$-vectors, or only all $d$-vectors, is injective and has rank $rn$. Indeed, a vector in such a partial kernel would have to be a linear combination of the displayed scaling directions; vanishing of the other two components and nonzeroness of the parameter vectors force all its coefficients to vanish.

Let $\mathcal U$ be representatives for $K/H$, including $0$. At the zero-extended point, the columns indexed by parameters in $H$ have support in $H\times H$ and reproduce the restricted Jacobian. For every $u\in\mathcal U\setminus\{0\}$ the remaining columns split into three groups with the following mutually disjoint supports:
\[
\begin{array}{c|c}
\text{varied parameters}&\text{output block}\\\hline
 a_{u+h},\ h\in H &(u+H)\times H\\
 b_{u+h},\ h\in H &H\times(u+H)\\
 d_{u+h},\ h\in H &(u+H)\times(u+H).
\end{array}
\]
No other output blocks occur in the differential. The same block supports hold simultaneously for every channel.

Each outside block is equivalent, by invertible row and column diagonal scalings, to the corresponding partial Jacobian on $H$. Here are those scalings explicitly. In the $a$-block write $i=u+i'$, $j=j'$, with $i',j'\in H$, and denote the outside parameter by $a_{u+h}$. The extra phase is
\[
 \beta_s(z,u)=\beta_s(i',u)\beta_s(h,u)^{-1},
 \qquad z=i'-h.
\]
Thus the row factor is $\beta_s(i',u)$ and the column factor is $\beta_s(h,u)^{-1}$. In the $b$-block the factors are $\beta_t(j',u)$ and $\beta_t(h,u)^{-1}$. In the $d$-block write $i=u+i'$, $j=u+j'$, $z=u+h$. The extra phase is
\[
 \beta_s(u,i'-h)\beta_t(u,j'-h)
 =\beta_s(u,i')\beta_t(u,j')
  [\beta_s(u,h)\beta_t(u,h)]^{-1}.
\]
These factors again separate into a row factor and a column factor, and all are nonzero. Therefore every outside partial block has rank $rn$.

All displayed blocks have disjoint row and column supports, so their ranks add. The complete rank is
\[
 r(3n-2)+3(|K/H|-1)rn=r(3N-2).
\]
The scaling argument of Lemma~\ref{lem:scaling-bound} also applies to $F_K$ and gives image dimension at most $r(3N-2)$. The computed Jacobian rank gives the reverse inequality.
\end{proof}

\begin{corollary}[One-channel dimension for finite abelian groups]\label{cor:abelian-one-channel}
Let $K$ be a finite abelian group and let $\beta_s,\beta_t:K\times K\to\C^\times$ be bicharacters. For the parameterization $F_K$ of Theorem~\ref{thm:subgroup-lift},
\[
\dim\overline{\im F_K}=3|K|-2.
\]
\end{corollary}
\begin{proof}
Take $H=\{0\}$ and $r=1$. Since $\beta_s(0,0)=\beta_t(0,0)=1$, the restricted map is $F_H(a,b,d)=abd$. Its differential at $(1,1,1)$ has rank $1=3|H|-2$. Theorem~\ref{thm:subgroup-lift} gives the assertion, including when $K$ is trivial.
\end{proof}

\begin{remark}[Necessity of the nonvanishing hypothesis]
The rank hypothesis in Theorem~\ref{thm:subgroup-lift} does not suffice if a parameter vector is zero. Indeed, take $H=\{0\}$, $r=1$, and $(a,b,d)=(0,1,1)$. The restricted differential has rank $1=3|H|-2$. After zero extension to a nontrivial group $K$, it is
\[
dF_K(x,y,z)_{i,j}=x_i\delta_{j,0}.
\]
Its rank is $|K|<3|K|-2$.
\end{remark}

\begin{remark}[Application to the twisted Heisenberg channel]
For $K=(\F_q^2,+)$, take
\[
 \beta_t((x,y),(u,v))=\psi_t(-xv).
\]
These are bicharacters, and the displayed polynomial channel is exactly the basepoint matrix parameterization $A(a)^T\operatorname{diag}(d)B(b)$. Hence Theorem~\ref{thm:subgroup-lift} applies to every central-character pair without an untwistedness assumption. The hypothesis requires rank $r(3|H|-2)$. A filling witness of smaller rank does not satisfy this hypothesis. Equality $r(3|H|-2)=|H|^2$ is allowed; the witnesses used below have rank strictly less than $|H|^2$.
\end{remark}

\begin{corollary}[Dimensions over extension fields]\label{cor:extension-families}
For every central-character pair, the expected affine secant dimension
\[
 \dim\overline{\operatorname{im}\Phi_r}=r(3q^2-2)
\]
holds in each of the following ranges:
\begin{enumerate}
\item $q=3^f$, $f\ge1$, and $1\le r\le3$;
\item $q=5^f$, $f\ge1$, and $1\le r\le8$;
\item $q=2^f$, with either $f=2$ or $f\ge4$, and $1\le r\le5$.
\end{enumerate}
For $q=3,4,5$ these are the certified cases of Theorem~\ref{thm:sharp-smallq}; for the remaining extension fields they are consequences of subgroup lifting.
\end{corollary}
\begin{proof}
For the first two assertions take $H=\F_p^2\subseteq\F_q^2$, where $p=3$ or $5$. Write the fixed nontrivial additive character as $\psi(z)=\exp(2\pi i\operatorname{Tr}_{\F_q/\F_p}(\kappa z)/p)$ for some $\kappa\ne0$. The restricted bicharacters are precisely the prime-field models with effective central types
\[
 s_0=\operatorname{Tr}_{\F_q/\F_p}(\kappa s),\qquad
 t_0=\operatorname{Tr}_{\F_q/\F_p}(\kappa t).
\]
Every such pair was certified at all sub-filling orders up to $3$ for $p=3$ and up to $8$ for $p=5$. A corresponding prefix of the frozen witness supplies the hypothesis of Theorem~\ref{thm:subgroup-lift}, by Lemma~\ref{lem:prefix-compression}. Lemma~\ref{lem:prefix-compression} also shows that each individual channel differential has rank $3|H|-2$. Each triple therefore has nonzero components: a zero component would reduce its differential rank to at most $|H|<3|H|-2$.

For characteristic two and $f\ge4$, choose a four-dimensional $\F_2$-linear subspace $U\subseteq\F_q$ and let $H=U\times\{0\}$. Then $H\cong (\F_2)^4\cong(\F_4^2,+)$ and both restricted bicharacters are trivial, for every ambient pair $(s,t)$. The untwisted $q=4$ certificates supply rank $r(3\cdot16-2)$ for $r\le5$. Apply the lifting theorem. The case $f=2$ is already part of the small-field theorem.
\end{proof}

\begin{remark}
Corollary~\ref{cor:extension-families} uses the exact certificates on the base groups and the subgroup-lifting theorem. Its ranges concern the stated secant orders and do not determine the first filling width for arbitrary $q$.
\end{remark}

\section{Secant geometry and comparison of decomposition widths}\label{sec:secant-geometry}

Let \(X_{s,t}\subset V_{s,t}\) be the affine cone obtained as the Zariski closure of one-channel maps.  We use the convention \(\sigma_r(\mathbb P X)\) for the Zariski closure of projective linear spans of \(r\) points of \(\mathbb P X\), so \(\sigma_1(\mathbb P X)=\mathbb P X\).  With this convention, the closure of the width-\(r\) image is the affine cone over \(\sigma_r(\mathbb P X_{s,t})\).  Thus the generic fusion width is the first secant order whose affine cone fills the ambient bilinear space.

By Theorem~\ref{thm:one-channel-dim},
\[
\dim \mathbb P X_{s,t}=3q^2-3
\]
for every type pair.  Hence the expected projective dimension of the \(r\)-th secant is
\begin{equation}
\operatorname{expdim}\sigma_r(\mathbb P X_{s,t})
=
\min\{r(3q^2-2)-1,\ q^4-1\}.
\label{eq:expected-secant-dim}
\end{equation}
Thus Conjecture~\ref{conj:general-fusion} concerns nondefectivity at the first expected filling order.

For the untwisted case, we obtain a toric description. Let
\[
A=G/Z\cong(\F_q^2,+),
\]
and take Fourier transform on the finite abelian group \(A\).  Ordinary convolution operators diagonalize.  If their Fourier multipliers are denoted by \(A_a,B_b,C_k\), then pointwise multiplication in the spatial domain becomes convolution in frequency.  Equivariance forces output frequency \(k=a+b\), and one channel has coefficient array
\begin{equation}
D_{a,b}=A_aB_bC_{a+b},
\qquad a,b\in\widehat A.
\label{eq:toric-tripod}
\end{equation}
After replacing the third index by its negative, \eqref{eq:toric-tripod} is the standard monomial parameterization of a group-based three-leaf claw (tripod) model for the finite abelian group \(A\).  Consequently:

\begin{proposition}[Untwisted toric correspondence]
The projective one-channel variety for \((s,t)=(0,0)\) is the group-based tripod toric variety for \(A\cong\F_q^2\), up to a relabeling of the third Fourier index.  The projectivizations of the width-\(r\) image closures are its \(r\)-th secant varieties.  The displayed Fourier formula depends on the normalization convention for the finite Fourier transform: with a normalized inverse transform a nonzero global scalar factor may appear in \(D_{a,b}\).  Such a factor is absorbed into the third kernel \(C\), and therefore does not alter the affine/projective one-channel variety or any fusion-width statement.
\end{proposition}

This identifies the untwisted fusion problem with a familiar nondefectivity question from algebraic statistics \citep{BanosEtAl2019}. More precisely, its all-orders nondefectivity is a special case of their Conjecture 5.1. Their Proposition 5.3 proves the \(q=2\) untwisted second-secant case, since \(A=\F_2^2\) gives the Kimura three-parameter tripod. We include that case for completeness and reproducibility, not as a new theorem.  If $s\ne0$ or $t\ne0$, at least one input commutant is isomorphic to the noncommutative algebra $M_q(\C)$. Thus the simultaneous scalar Fourier diagonalization used for $(s,t)=(0,0)$ does not apply to these pairs.  Theorem~\ref{thm:second-secant} proves order-two nondefectivity uniformly in these twisted sectors. The exact certificates reach the first filling order for every prime power $q\le5$, and Theorem~\ref{thm:subgroup-lift} propagates their sub-filling ranks to the infinite families above.

For the polynomial parameterization used here, filling of the secant closure is equivalent to dominance. In either formulation, constructibility supplies a nonempty Zariski-open subset of the actual image, so the generic statement is about exact width. Exact width and border width may still differ at special targets. The full-rank Jacobian certificates provide an explicit sufficient test for this dominance.

\begin{remark}[Comparison with quadratic Waring rank]\label{rem:waring-comparison}
The two notions of width differ even for one fusion channel. For any type pair, take $A=B=I$ and choose the output kernel so that its basepoint row is $d=\mathbf 1$. Equation~\eqref{eq:channel-matrix} then gives $M=I_N$, where $N=q^2$. This nonzero equivariant bilinear map has fusion width one. Its scalar basepoint polynomial on $\mathcal H_s\oplus\mathcal H_t\cong\C^{2N}$ is
\[
 p(f,g)=f^Tg
 =\begin{pmatrix}f\\g\end{pmatrix}^{\!T}
 \begin{pmatrix}0&I_N/2\\I_N/2&0\end{pmatrix}
 \begin{pmatrix}f\\g\end{pmatrix}.
\]
The symmetric coefficient matrix has rank $2N$, so at least $2N$ squares of linear forms are required. The identity
\[
 p(f,g)=\frac14\sum_{i=1}^{N}\bigl((f_i+g_i)^2-(f_i-g_i)^2\bigr)
\]
attains this bound over $\C$. Thus the quadratic Waring rank is $2N$, as in the symmetric-rank characterization of Akbari and Jamshidi (2026, Theorem~6). Here the Waring rank refers to the scalar basepoint polynomial, whereas fusion width refers to its equivariant vector-valued extension.
\end{remark}

\section{A conjecture on generic width}\label{sec:width-conjecture}

Theorem~\ref{thm:sharp-smallq} leads to the following conjecture.

\begin{conjecture}[Generic fusion width]\label{conj:general-fusion}
For every prime power \(q\) and every \(s,t\in\F_q\),
\begin{equation}
 r_{\rm gen}(q;s,t)
 =
 \left\lceil\frac{q^4}{3q^2-2}\right\rceil.
\label{eq:main-conjecture}
\end{equation}
\end{conjecture}

For $q=2,3,4,5$, Theorem~\ref{thm:sharp-smallq} proves the conjecture and also gives the expected dimensions at every smaller secant order. These cases include the extension field $\F_4$. They do not exclude the possibility of a defect at a larger prime power.

By Theorem~\ref{thm:one-channel-dim}, Conjecture~\ref{conj:general-fusion} is equivalent to nondefectivity at the first expected filling order. It makes no assertion about every smaller secant order. In the untwisted case, the question belongs to the theory of group-based toric mixtures. For nonzero central characters, the commutants are matrix algebras, so the same Fourier reduction is unavailable. Theorem~\ref{thm:second-secant} settles order two uniformly, and Corollary~\ref{cor:extension-families} settles further orders on infinite families. Reaching the first filling order for arbitrary $q$ requires additional control of sums of more differential images.

The lower bound satisfies
\[
\frac{q^4}{3q^2-2}=\frac{q^2}{3}+O(1).
\]
Thus generic representation requires a number of channels that grows quadratically with $q$. The width is measured with respect to the prescribed channel map. The certificates do not give an algorithm attaining the generic minimum for an arbitrary input, nor do they establish accuracy or efficiency on a learning task.

\section{Reproducibility}\label{sec:reproducibility}

Online Resource~1 contains stored parameter points, specified nonzero minors and certificates for all secant dimensions, and read-only verification scripts. Table~\ref{tab:reproduction} summarizes the entry points. The verifiers compare the computed values with the stored indices and determinants and do not modify the certificate files.
\begin{table}[ht]
\caption{Verification and constructive-reconstruction entry points in Online Resource~1.}
\label{tab:reproduction}
\centering\small
\begin{tabular}{p{0.46\textwidth}p{0.45\textwidth}}
\toprule
Entry point & Purpose\\\midrule
\nolinkurl{verify_all.py} & Check hashes, run the certificate suite and confirm unchanged file bytes\\
\nolinkurl{verify_checksums.py} & Validate the SHA-256 manifest before mathematical verification\\
\nolinkurl{verify_q2_all_twisted_sectors.py}, \nolinkurl{verify_q3_all_twisted_sectors.py} & Verify all prefix ranks and fixed terminal determinants\\
\nolinkurl{verify_q45_frozen_certificates.py} & Verify stored hashes, penultimate ranks and stored terminal minors\\
\nolinkurl{verify_full_secant_profiles_q45.py} & Recompute every stored sub-filling minor for \(q=4,5\)\\
\nolinkurl{verify_independent_coordinates.py} & Reconstruct the Jacobians directly from coefficient sums, independently of the matrix-product implementation\\
\nolinkurl{decompose_basepoint_matrix.py} & Exact rational rank-factorization example for Proposition~\ref{prop:constructive-upper}\\
\nolinkurl{verify_theory_upgrades.py} & Exact checks of the second-secant parameter points and a subgroup lift\\
\bottomrule
\end{tabular}

\end{table}

The tested environment was Python~3.12.14 and NumPy~2.3.5 on Windows 10. The \(q=2,3\) verifiers use Python integer arithmetic. The \(q=4,5\) verifiers use integer arrays reduced modulo \(101\), with modular inverses; no floating-point rank routines enter the certificates. For these fixed dimensions, unreduced dot products have magnitude at most \(25\cdot100^3\), well within the signed 64-bit range. Integer parameter residues are lifted to ordinary integers and character coefficients via the displayed cyclotomic specialization; no embedding of a finite field into \(\C\) is assumed.

To reproduce the computations, run \nolinkurl{python verify_all.py}. It validates checksums before and after the verification runs, refuses optimized Python execution that would disable legacy assertions, and reports nonzero exit status on a failed check. The independent coordinate verifier is an optional additional audit. The script \nolinkurl{verify_theory_upgrades.py} checks 32 stored two-channel parameter points for $q=2,3,4,5,7,8,9,11$, compares two Jacobian constructions, and verifies four width-three lifts from $q=3$ to $q=9$, each of rank $723$. The proofs of Theorems~\ref{thm:second-secant} and~\ref{thm:subgroup-lift} do not depend on these checks. They use Python-integer modular elimination over $\F_{2311}$; small matrix-assembly comparisons at $q=3,4$ use NumPy integer arrays within the signed 64-bit range. The witness SHA-256 values are
\begin{quote}\scriptsize\ttfamily
q4: affb0a1a77348bb33f7ac797529fa152958f8cdbe19616740e8dbf59811e3c15\\
q5: 3a6d802a5b097dd5d8081a17c2dba3303ba7792156b06236e87f2c74557d4e19.
\end{quote}
All certificate data are self-contained and require no external datasets. The full profile can be proved using only the penultimate and terminal rank checks by Lemma~\ref{lem:prefix-compression}; the additional stored minors make each order separately inspectable.

\section{Conclusion and future work}\label{sec:conclusion}

We have determined the one-channel dimension and the second-secant dimension for the bilinear channel varieties associated with representations of $UT_3(\F_q)$ induced from central characters. The affine one-channel variety has dimension $3q^2-2$ for every type pair. The affine second-secant closure has dimension $6q^2-4$ for $q\ge3$, while it fills the ambient space for $q=2$. The subgroup-lifting theorem transfers differential rank $r(3|H|-2)$ on $H$ to rank $r(3|K|-2)$ on $K$, provided that every parameter vector is nonzero. Together with the exact small-field certificates, it gives higher-secant dimensions on the extension-field families in Corollary~\ref{cor:extension-families}. For $q=2,3,4,5$, all secant dimensions are determined and the generic widths are $2,4,6,9$, respectively.

The main unresolved question is Conjecture~\ref{conj:general-fusion}, which asserts that the dimension lower bound $\lceil q^4/(3q^2-2)\rceil$ is attained for every prime power and every central-character pair. A proof requires control of the differential rank at the first expected filling order. The present subgroup construction preserves a fixed channel order and therefore does not reach that order when the extension degree grows. One possible approach is to combine several translated subgroup supports and determine the intersections of their differential images. This would require an argument beyond the zero-extension construction used here. Determining the dimensions at every smaller secant order is a separate, stronger problem.

A further question concerns reconstruction at small width. The rank-factorization construction gives at most $q^2$ channels for any target, whereas the generic-width results establish existence without giving a procedure that attains the generic minimum. An algorithm for finding such decompositions would require analysis of solvability and numerical stability. The nonsingular minor used in the one-channel proof gives local information at the selected parameter point; it does not provide a conditioning estimate for general decompositions. The polynomial realization and approximation results in Appendices~\ref{app:polynomial} and~\ref{app:approximation} do not supply such quantitative decomposition bounds.

\section*{Supplementary information}
Online Resource~1 (\nolinkurl{ESM_1.zip}): exact computational certificates and verification code for the secant dimensions, with stored parameter points, nonzero-minor data, reproduction logs and a SHA-256 checksum manifest.

\appendix
\section{Details of the basepoint matrix formula}

Let \(e_i\), \(1\le i\le N\), be the delta basis of \(\C^{\F_q^2}\).  Write
\[
Af=A f,\qquad Bg=B g
\]
in matrix notation, and let \(c_o\) be the basepoint row of \(C\).  Then
\[
\Gamma_{A,B,C}(f,g)(o)
=
\sum_{z=1}^N (c_o)_z (Af)_z(Bg)_z.
\]
Expanding gives
\[
\sum_{i,j}f_i g_j
\left(
\sum_z A_{z i}(c_o)_zB_{zj}
\right),
\]
so the coefficient matrix is
\[
M=A^T\operatorname{diag}(c_o)B,
\]
which is \eqref{eq:channel-matrix}.  Differentiating this expression with respect to the twisted-shift kernel coefficients gives the Jacobian implemented in the supplementary scripts.

\section{Dominance and exact representation}

For completeness, let \(X=\A^m\), \(Y=\A^n\), and \(\Phi:X\to Y\) be a polynomial map over \(\C\).  If \(d\Phi_x\) has rank \(n\) at one point, then the image has dimension \(n\), hence \(\Phi\) is dominant.  Chevalley's theorem says \(\Phi(X)\) is constructible.  A dense constructible subset of irreducible \(Y\) contains a nonempty open subset.  Thus dominance obtained from a Jacobian certificate implies that a generic target is represented by an actual finite parameter point; no limiting argument is needed.

\section{Certificate orbit counts}

For reference, simultaneous scaling and exchange classify nonzero ratios by
\[
r\sim r^{-1},\qquad r=t/s\in\F_q^\times.
\]
Pairs with exactly one zero form one orbit, and \((0,0)\) forms another.  Thus the representative lists used in Section~\ref{sec:type-symmetry} are exhaustive for \(q=2,3,4,5\).

\section{Equivariant polynomial maps}\label{app:polynomial}

The neural-network interpretation uses group-equivariant convolutions as in Cohen and Welling~\citep{CohenWelling2016} and Kondor and Trivedi~\citep{KondorTrivedi2018}. General equivariant universality results are available in Ravanbakhsh~\citep{Ravanbakhsh2020} and Nyholm et al.~\citep{NyholmEtAl2026}. In \emph{Graded Representation Theory of Equivariant Neural Networks}, Shayestehfar (2026; \href{https://arxiv.org/abs/2609.25776}{arXiv:2609.25776}, Theorem~3.6 and Proposition~5.1) describes the degree components and their span for fixed equivariant lifts and readouts. His Proposition~5.3 proves that regular-representation branches span every finite-group real equivariant homogeneous polynomial space. Appendix~\ref{app:polynomial} gives a direct implementation of this existence principle by central-type linear maps, products and conjugation. The width results in Sections~5--9 concern the dimensions obtained when all channel parameters vary.

We now apply the preceding description of the commutants to equivariant polynomial maps. In this section, depth and width are unrestricted.

Let
\[
E=\bigoplus_{j=1}^m \mathcal H_{t_j}^{\oplus d_j},
\qquad
F=\bigoplus_{k=1}^n \mathcal H_{u_k}^{\oplus e_k}.
\]
A $*$-polynomial is a polynomial in the complex coordinates and their conjugates. A \emph{typed polynomial network} is built from the following operations:
\begin{enumerate}[label=(\roman*)]
\item arbitrary \(G\)-equivariant complex-linear maps between channels of the same type;
\item pointwise products \(\mathcal H_s\times\mathcal H_t\to\mathcal H_{s+t}\);
\item complex conjugation \(\mathcal H_t\to\mathcal H_{-t}\);
\item additions, scalar multiplication, and spatially constant fields of type \(0\).
\end{enumerate}
A constant node is a scalar multiple of the all-ones field in \(\mathcal H_0\); an arbitrary fixed nonconstant field would not be equivariant. Conjugation is included only for \(*\)-polynomial and real-continuous approximation; purely holomorphic polynomial maps do not require it.
Formally, such a network is a finite directed acyclic computation graph.  Every edge carries a finite direct sum of typed feature channels; linear nodes are intertwiners, product nodes have two incoming typed channels and one outgoing channel of summed type, and conjugation nodes reverse the type.  Fan-out is allowed, so a previously computed feature may feed several later nodes.  No limiting operation, division, or data-dependent branching is permitted.  Exact representation below means equality with the map computed by such a graph.

\begin{lemma}[Equivariant maps are determined at one basepoint]
Let \(Q:E\to\mathcal H_u\) be \(G\)-equivariant.  If \(Q(v)(o)=0\) for every \(v\in E\), then \(Q=0\).
\end{lemma}

\begin{proof}
Every point of \(G/Z\) is \(g\cdot o\) for some \(g\in G\).  Equivariance expresses \(Q(v)(g\cdot o)\) in terms of the basepoint value of \(Q(g^{-1}v)\), up to the known fiber action.  Hence all coordinates vanish.
\end{proof}

The general spanning result of Shayestehfar (2026, Proposition~5.3) also yields existence in the present network class. Apply it degree by degree to the underlying real representations and then complexify. The complexified input is $E\oplus\overline E$, with the original input represented by $(v,\overline v)$. The complexified output is $F\oplus\overline F$, and projection onto $F$ recovers the original map. The regular hidden space decomposes as $\C[G]\cong\bigoplus_{t\in\F_q}\mathcal H_t$. Its equivariant linear maps split by central type, and its pointwise powers expand into products $\mathcal H_s\mathcal H_t\subseteq\mathcal H_{s+t}$. Conjugation supplies the conjugate input factors, and constants have type zero. Thus these branches can be implemented by the allowed operations. The proof below gives a direct basepoint construction. Neither argument gives a minimal-width decomposition.

\begin{theorem}[Exact realization of equivariant \(*\)-polynomials]
Every \(G\)-equivariant \(*\)-polynomial map \(Q:E\to F\) is represented exactly by a finite typed polynomial network.
\label{thm:poly-realization}
\end{theorem}

\begin{proof}
It is enough to consider one output copy \(\mathcal H_u\).  Set
\[
p(v)=Q(v)(o).
\]
This is a scalar \(*\)-polynomial in the input coordinates.  Because \(z\in Z\) fixes the basepoint and acts on the output by \(\psi_u(z)\), equivariance gives
\[
p(z\cdot v)=\psi_u(z)p(v).
\]
Thus only central-weight-\(u\) monomials occur in the character decomposition of \(p\).

A constant term can occur only when $u=0$ and is represented by a constant node. Consider a monomial of positive degree. Each holomorphic factor is an input coordinate
\(f(x,y)\) of some type \(t\), and each antiholomorphic factor is its conjugate.  By \eqref{eq:twisted-shift},
\[
(R_{x,y}^tf)(o)=f(x,y).
\]
Therefore every holomorphic coordinate factor is obtained at the basepoint by an allowed typed linear operator; an antiholomorphic factor is obtained by following the shift with conjugation, which changes type \(t\) to \(-t\).  Repeated pointwise multiplication realizes the entire monomial as an \(\mathcal H_u\)-valued equivariant feature because its total central weight is \(u\).  Summing the monomial subnetworks realizes a map whose basepoint value is \(p\).  By the previous lemma it equals \(Q\) everywhere.  Repeating this for each output channel proves the theorem.
\end{proof}

\begin{remark}
In this construction, the type of each monomial is the sum of the types of its holomorphic factors minus the sum of the types of its antiholomorphic factors.
\end{remark}

\section{Universal approximation}\label{app:approximation}

The polynomial-realization theorem gives the following approximation result by Stone--Weierstrass approximation and Reynolds averaging.

Equip all finite-dimensional representations with \(G\)-invariant Hermitian norms.  Let \(K\subset E\) be compact and \(G\)-invariant.

\begin{theorem}[Universal approximation]
For every continuous \(G\)-equivariant map \(F_0:K\to F\) and every \(\varepsilon>0\), there exists a typed polynomial network \(N: E\to F\) such that
\[
\sup_{v\in K}\|N(v)-F_0(v)\|<\varepsilon.
\]
\end{theorem}

\begin{proof}
The algebra generated by complex coordinate functions and their conjugates is self-adjoint, contains constants, and separates points of the compact subset \(K\).  By the complex Stone--Weierstrass theorem, vector-valued \(*\)-polynomials approximate \(F_0\) uniformly coordinatewise.  Choose such a polynomial map \(P:E\to F\) with uniform error less than \(\varepsilon\) on \(K\).

Apply Reynolds averaging:
\begin{equation}
\mathcal R P(v)
=
\frac1{|G|}\sum_{g\in G}
\rho_F(g)P(\rho_E(g^{-1})v).
\label{eq:reynolds}
\end{equation}
This is again a \(*\)-polynomial and is exactly \(G\)-equivariant.  Since \(K\) is invariant and the representations are unitary,
\[
\|\mathcal RP(v)-F_0(v)\|
\le
\frac1{|G|}\sum_g
\|P(\rho_E(g^{-1})v)-F_0(\rho_E(g^{-1})v)\|
<\varepsilon.
\]
Theorem \ref{thm:poly-realization} realizes \(\mathcal RP\) exactly by a typed polynomial network.
\end{proof}

The argument combines Theorem~\ref{thm:poly-realization} with polynomial approximation and averaging. General equivariant approximation results are given in Ravanbakhsh~\citep{Ravanbakhsh2020} and Nyholm et al.~\citep{NyholmEtAl2026}.

\section*{Statements and Declarations}

\subsection*{Funding}
No funding was received for conducting this study.

\subsection*{Author contributions}
Nguyen Bich Van developed the theoretical foundations and proofs and wrote the verification code. Nguyen Trung Kien carried out computations for specific examples and checked the code.

\subsection*{Competing interests}
The authors declare that they have no competing financial or non-financial interests that are directly or indirectly related to this work.

\subsection*{Data and code availability}
The computational data and code supporting the results are supplied as Online Resource~1. They include the exact computer-assisted certificates and frozen witnesses used for Theorem~\ref{thm:sharp-smallq}, verification scripts, reproduction logs, and a SHA-256 checksum manifest. All computations use exact arithmetic; no external datasets were used.

\begingroup
\interlinepenalty=10000
\bibliography{references}
\endgroup

\end{document}